\documentclass[12pt]{amsart}

\usepackage[latin1]{inputenc}
\usepackage[T1]{fontenc}
\usepackage{amsfonts,amssymb,latexsym}
\usepackage{amsthm}
\usepackage[dvips]{color}

\usepackage{url}

\usepackage{float}
\usepackage{graphicx}

\newtheorem{theorem}{Theorem}[section]
\newtheorem{lemma}[theorem]{Lemma}
\newtheorem{prop}[theorem]{Proposition}
\newtheorem{cor}[theorem]{Corollary}

\theoremstyle{definition}
\newtheorem{defi}[theorem]{Definition}

\theoremstyle{remark}
\newtheorem{remark}[theorem]{Remark}

\numberwithin{equation}{section}

\newcommand{\rr}{{\mathbb R}}

\newcommand{\rd}{{\mathbb R^d}}

\newcommand{\supp}{\operatorname{supp}}
\newcommand{\spann}{\operatorname{span}}
\newcommand{\ip}[2]{{\langle #1,#2\rangle}}
\renewcommand{\Re}{\operatorname{Re}}

\newcommand{\eqfd}{\stackrel{f.d.}{=}}
\newcommand{\eqd}{\stackrel{d}{=}}
\newcommand{\Exp}{\mathbb E}

\newcommand{\diag}{\operatorname{diag}}

\newcommand{\trace}{\operatorname{trace}}

\begin{document}
\sloppy
\title[Operator Scaling Stable Fields]{Operator Scaling Stable Random Fields}

\author{Hermine Bierm\'e}\thanks{Hermine Bierm\'e was supported by NSF grant DMS-0417869.}
\address{Hermine Bierm\'e, MAP5 Universit\'e Ren\'e Descartes, 45 rue des Saints Pères, 75270 Paris cedex 06, France}
\email{hermine.bierme\@@{}math-info.univ-paris5.fr}
\author{Mark M. Meerschaert}\thanks{Mark M. Meerschaert was partially supported by NSF grants DMS-0417869 and DMS-0139927, and Marsden grant UoO-123 from the Royal Society of New Zealand.}
\address{Mark M. Meerschaert, Department of Mathematics \& Statistics,
 University of Otago, New Zealand}
\email{mcubed\@@{}maths.otago.ac.nz}
\urladdr{http://www.maths.otago.ac.nz/~mcubed/}
\author{Hans-Peter Scheffler}\thanks{Hans-Peter Scheffler was partially supported by NSF grant DMS-0417869.}
\address{Hans-Peter Scheffler, Fachbereich Mathematik, Universit\"at Siegen, 57068 Siegen, Germany}
\email{peter.scheffler\@@{}mac.com}

\date{20 February 2006}

\begin{abstract}
A scalar valued random field $\{X(x)\}_{x\in\rd}$ is called {\it
operator-scaling} if for some $d\times d$ matrix $E$ with positive
real parts of the eigenvalues and some $H>0$ we have
\begin{equation*}
\{X(c^Ex)\}_{x\in\rd}\eqfd\{c^HX(x)\}_{x\in\rd}\quad\text{for all
$c>0$,}
\end{equation*}
where $\eqfd$ denotes equality of all finite-dimensional marginal
distributions. We present a moving average and a harmonizable
representation of stable operator scaling random fields by
utilizing so called $E$-homogeneous functions $\varphi$,
satisfying $\varphi(c^Ex)=c\varphi(x)$. These fields also have
stationary increments and are stochastically continuous. In the
Gaussian case critical H\"older-exponents and the
Hausdorff-dimension of the sample paths are also obtained.
\end{abstract}


\maketitle

\baselineskip=18pt

\section{Introduction}

A scalar valued random field $\{X(x)\}_{x\in\rd}$ is called {\it
operator-scaling} if for some $d\times d$ matrix $E$ with positive
real parts of the eigenvalues and some $H>0$ we have
\begin{equation}\label{eq11}
\{X(c^Ex)\}_{x\in\rd}\eqfd\{c^HX(x)\}_{x\in\rd}\quad\text{for all
$c>0$,}
\end{equation}
where $\eqfd$ denotes equality of all finite-dimensional marginal
distributions. As usual $c^E=\exp(E\log c)$ where
$\exp(A)=\sum_{k=0}^\infty \frac{A^k}{k!}$ is the matrix exponential. Note
that if $E=I$, the identity matrix, then \eqref{eq11} is just the
well-known {\it self-similarity} property
$\{X(cx)\}_{x\in\rd}\eqfd\{c^HX(x)\}_{x\in\rd}$ where one usually
calls $H$ the Hurst-index. See \cite{em:mae} for an overview of
self-similar processes in the one-dimensional case $d=1$.
Self-similar processes are used in various fields of applications
such as internet traffic modelling \cite{wil:taqqu}, ground water modelling and
mathematical finance, just to mention a few. Various examples can be found for instance in the books \cite{levyvehel} and  \cite{Abry}. A very important
class of such fields or processes are the {\it fractional stable
fields} and especially the L\'evy fractional Brownian field.

These fields have different definitions which are usually not
equivalent. More precisely, for $0<\alpha\leq 2$ let
$Z_\alpha(dy)$ be an independently scattered symmetric $\alpha$-stable ($S\alpha S$) random
measure on $\rd$ with Lebesgue control measure $\lambda^d$  (see \cite{sam:taqqu} p.\ 121). For
$0<H<1$ one defines the {\it moving average} representation by
\begin{equation}\label{eq12}
X_H(x)=\int_\rd\Bigl(\|x-y\|^{H-d/\alpha}-\|y\|^{H-d/\alpha}\Bigr)\,Z_\alpha(dy).
\end{equation}
For $W_{\alpha}(d\xi)$ a complex isotropic $S\alpha S$ random
measure with Lebesgue control measure the harmonizable
representation is given by
\begin{equation}\label{eq13}
\tilde X_H(x)=\Re\int_\rd\bigl(e^{i\ip
x\xi}-1\bigr)\|\xi\|^{-H-d/\alpha}\,W_{\alpha}(d\xi) .
\end{equation}
See \cite{sam:taqqu} for a comprehensive introduction to random
integrals with respect to stable measures. It follows from basic
properties that $\{X_H(cx)\}_{x\in\rd}\eqfd\{c^HX(x)\}_{x\in\rd}$
as well as $\{\tilde X_H(cx)\}_{x\in\rd}\eqfd\{c^H\tilde
X_H(x)\}_{x\in\rd}$. Moreover, both processes have {\it stationary
increments}, that is for any $h\in\rd$ we have
$\{X_H(x+h)-X_H(h)\}_{x\in\rd}\eqfd\{X_H(x)\}_{x\in\rd}$ and
similarly for $\{\tilde X_H(x)\}_{x\in\rd}$. Furthermore both
fields are {\it isotropic}, that is
$\{X_H(Ax)\}_{x\in\rd}\eqfd\{X(x)\}_{x\in\rd}$ for any orthogonal
matrix $A$. It is worth mentioning that if $\alpha<2$ the fields
$\{X_H(x)\}_{x\in\rd}$ and $\{\tilde X_H(x)\}_{x\in\rd}$ defined
in \eqref{eq12} and \eqref{eq13}, respectively, are usually
different. See \cite{sam:taqqu}, Theorem 7.7.4 for the
one-dimensional case. However, in the Gaussian case $\alpha=2$, by
computing the covariance function of the fields, it follows that
$\{X_H(x)\}_{x\in\rd}$ and $\{\tilde X_H(x)\}_{x\in\rd}$ have the same law
up to a multiplicative constant and known as the L\'evy fractional
Brownian field.

Certain applications (see, e.g., \cite{benson,Bonami,Ponson} and references therein)
require that the random field is anisotropic and satisfies a
scaling relation. This scaling relation should have different
Hurst indices in different directions and these directions should
not necessarily be orthogonal. In the Gaussian case a prominent
example of an anisotropic random field is the fractional Brownian
sheet $\{B_H(x)\}_{x\in\rd}$ defined as follows: Let $0<H_j<1$ for
$j=1,\dots,d$ and set $H=(H_1,\dots,H_d)$. Define
\[
B_H(x)=\int_\rd\prod_{j=1}^d\bigl[|x_j-u_j|^{H_j-1/2}-|u_j|^{H_j-1/2}\bigr]\,Z_2(du)
.\] See \cite{Xiao,herbin,xioa:zhang} and the literature cited
there for more information on these fields. Then, if we set
$E=\diag(H_1^{-1},\dots,H_d^{-1})$, it follows by a simple
computation that
$\{B_H(c^Ex)\}_{x\in\rd}\eqfd\{cB_H(x)\}_{x\in\rd}$ so
$\{B_H(x)\}_{x\in\rd}$ is operator scaling in the sense of
\eqref{eq11}. However, $\{B_H(x)\}_{x\in\rd}$ does not have
stationary increments.

The purpose of this paper is to define two different classes of
{\it operator scaling stable random fields} (OSSRF) and analyze
their basic properties. We present a moving average representation
as well as a harmonizable representation. Our constructions are
based on a class of $E$-homogeneous functions $\varphi :
\rd\to[0,\infty)$ where $\varphi$ is positive on
$\rd\setminus\{0\}$ and $\varphi(c^Ex)=c\varphi(x)$ for all
$x\in\rd$ and $c>0$. Such functions were studied in detail in
\cite{thebook}, Chapter 5. It will turn out that the harmonizable
representation allows more flexibility in the class of possible
functions $\varphi$ in contrast to the moving average
representation which is more restrictive. However, in both cases
the OSSRFs satisfy \eqref{eq11}, have stationary increments and
are continuous in probability. In the Gaussian case $\alpha=2$ we
show that there exists modifications of these fields which are
almost surely H\"older-continuous of certain indices and we
compute the Hausdorff-dimension of the graph.

This paper is organized as follows: In section 2 we introduce the
class of $E$-homogeneous functions, derive some basic properties
and provide important examples. In section 3 we define and analyze
a moving average representation of OSSRFs. Section 4 is devoted to
the harmonizable representation and its properties. Finally, in
the Gaussian case $\alpha=2$, we analyze the sample path
properties of both the moving average and the harmonizable
representation of OSSRFs.

\section{$E$-homogeneous functions}\label{Ehom}
Let $E$ be a real $d\times d$ matrix with positive
real parts of the eigenvalues $0<a_1<a_2\ldots <a_p$ for $p\le d$.
Let us denote $\Gamma=\rd\setminus\{0\}$. It follows from Lemma 6.1.5 of \cite{thebook} that there exists a
norm $\|\cdot\|_0$ on $\rd$ such that for the unit sphere
$S_0=\{x\in\rd : \|x\|_0=1\}$ the mapping $\Psi : (0,\infty)\times
S_0\to\Gamma$, $\Psi(r,\theta)=r^E\theta$ is a homeomorphism.
Moreover for any $x\in\Gamma$ the function $t\mapsto \|t^Ex\|_0$
is strictly increasing. Hence we can write any $x\in\Gamma$
uniquely as $x=\tau(x)^El(x)$ for some {\it radial part}
$\tau(x)>0$ and some direction $l(x)\in S_0$ such that
$x\mapsto\tau(x)$ and $x\mapsto l(x)$ are continuous. Observe that
$S_0=\{x\in\rd : \tau(x)=1\}$ is compact. Moreover we know that
$\tau(x)\to\infty$ as $x\to\infty$ and $\tau(x)\to 0$ as $x\to 0$.
Hence we can extend $\tau(\cdot)$ continuously by setting
$\tau(0)=0$. Note that further $\tau(-x)=\tau(x)$ and
$l(-x)=-l(x)$. The following result gives bounds on the growth
rate of $\tau(x)$ in terms of the real parts of the eigenvalues of
$E$.

\begin{lemma}\label{lem10}
For any (small) $\delta>0$ there exist constants $C_1,\dots,C_4>0$
such that
for all $\|x\|_0\leq 1$  or all $\tau(x)\leq 1$,
$$C_1\|x\|_0^{1/a_1+\delta}\le \tau(x) \leq C_2\|x\|_0^{1/a_p-\delta},$$
and, for all $\|x\|_0\geq 1$  or all $\tau(x)\geq 1$,
$$C_3\|x\|_0^{1/a_p-\delta} \le \tau(x) \leq C_4\|x\|_0^{1/a_1+\delta}.$$
\end{lemma}

\begin{proof}
We will only prove the first two inequalities. It follows from
Theorem 2.2.4 of \cite{thebook} that for any $\delta'>0$ we have
$t^{a_1-\delta'}\|t^{-E}\theta\|_0\to 0$ as $t\to\infty$ uniformly
in $\|\theta\|_0=1$. Hence $\|t^{-E}\|_0:=\sup_{\theta\in
S_0}\|t^{-E}\theta\|_0\leq Ct^{-a_1+\delta'}$ for all $t\geq 1$
and some constant $C>0$. Equivalently $\|s^E\|_0\leq
Cs^{a_1-\delta'}$ for all $s\leq 1$. Since
$\|x\|_0=\bigl\|\tau(x)^El(x)\bigr\|_0\leq \|\tau(x)^E\|_0\leq
C\tau(x)^{a_1-\delta'}$ we get $\tau(x)\geq
C_1\|x\|_0^{1/a_1+\delta}$, for
$\delta=\frac{1}{a_1-\delta_1}-\frac{1}{a_1}$, if $\|x\|_0\leq 1$
which is equivalent
to $\tau(x)\leq 1$.\\
 Similarly we know that, for any $\delta'>0$,
$t^{-a_p-\delta'}\|t^E\theta\|_0\to 0$ as $t\to\infty$ uniformly
in $\|\theta\|_0=1$. Therefore $\|t^E\|_0\leq Ct^{a_p+\delta'}$
for all $t\geq 1$ or equivalently $\|s^{-E}\|_0\leq
Cs^{-a_p-\delta'}$ for all $s\leq 1$. But
$x=\tau(x)^El(x)$ and $l(x)=\tau(x)^{-E}x$. Thus, $1\le \bigl\|\tau(x)^{-E}\bigr\|_0\bigl\|x\bigr\|_0$ and $\bigl\|x\bigr\|_0\geq
C^{-1}\tau(x)^{a_p+\delta'}$ for all $\|x\|_0\leq 1$. Hence
$\tau(x)\leq C_2\|x\|_0^{1/a_p-\delta}$ for
$\delta=\frac{1}{a_p}-\frac{1}{a_p+\delta'}$ and $\|x\|_0\leq 1$.
The proof is complete.
\end{proof}

The following results generalize some of the results in
\cite{folland}, Chapter 1.A to our more general case of exponents
$E$.

\begin{lemma}\label{lem11}
There exists a constant $K\geq 1$ such that for all $x,y\in\rd$ we
have
\[ \tau(x+y)\leq K\bigl(\tau(x)+\tau(y)\bigr) .\]
\end{lemma}

\begin{proof}
Observe that the set $G=\{(x,y)\in\rd\times\rd :
\tau(x)+\tau(y)=1\}$ is bounded by Lemma \ref{lem10} and closed by continuity of $\tau$.  Hence $G$ is a compact set.
Thus the continuous function $(x,y)\mapsto\tau(x+y)$ assumes a
finite maximum $K$ on $G$. Since $S_0\times\{0\}\subset G$, we have
$K\geq 1$. Given any $x,y\in\rd$ both not equal to zero we set
$s=\left(\tau(x)+\tau(y)\right)^{-1}$. Then, with $\tau(c^Ex)=c\tau(x)$ it follows
that \[
\tau(x+y)=s^{-1}\tau\bigl(s^E(x+y)\Bigr)=s^{-1}\tau\bigl((s^Ex)+(s^Ey)\bigr).\]
But $\left(s^Ex,s^Ey\right)\in G$ since $\tau\left(s^Ex\right)+\tau\left(s^Ey\right)=s\left(\tau(x)+\tau(y)\right)=1$.
Therefore,
$$\tau(x+y) \leq Ks^{-1}=K\bigl(\tau(x)+\tau(y)\bigr)$$ and the proof is complete.
\end{proof}

Now let $q=\trace(E)$ and observe that by multivariable change of
variables we have $\lambda^d(c^E(B))=c^q\lambda^d(B)$ for all
Borel sets $B\subset\rd$, $c>0$, which can be written as $d(c^Ex)=c^qdx$. Let
$B(r,x)=\{y\in\rd : \tau(y-x)<r\}$ denote the {\it ball} of radius
$r>0$ around $x\in\rd$. Then it is easy to see that
$B(r,x)=x+B(r,0)=x+r^EB(1,0)$ and hence
$\lambda^d(B(r,x))=r^q\lambda^d(B(1,0))$. The following
proposition provides an integration in {\it polar coordinates}
formula.

\begin{prop}\label{prop12}
There exists a unique finite Radon measure $\sigma$ on $S_0$ such
that for all $f\in L^1(\rd,dx)$ we have
\[\int_\rd f(x)\,dx=\int_0^\infty\int_{S_0}f(r^E\theta)\,\sigma(d\theta)\,r^{q-1}\,dr
.\]
\end{prop}
The proof of Proposition \ref{prop12} is based on the following.
\begin{lemma}\label{lem13}
If $f:\Gamma\to\mathbb C$ is continuous and $f(r^Ex)=r^{-q}f(x)$
for all $r>0$ and $x\in\Gamma$, then there exists a constant
$\mu_f$ such that for all $g\in
L^1\bigl((0,\infty),r^{-1}dr\bigr)$ we have
\[ \int_\rd f(x)g(\tau(x))\,dx=\mu_f\int_0^\infty g(r)\frac{dr}r
.\]
\end{lemma}

\begin{proof}
Let $L_f : (0,\infty)\to\mathbb C$ be defined as
\[ L_f(r)=\begin{cases} \int_{1\leq\tau(x)\leq r}f(x)\,dx &
\text{if $r\geq 1$} \\ -\int_{1\leq\tau(x)\leq r^{-1}}f(x)\,dx &
\text{if $r<1$ .} \end{cases} \] Since $f$ is continuous on $\Gamma$, from
dominated convergence, $L_f$ is continuous on $(0,1)\cup(1,+\infty)$.
But $\lambda^d(B(r,0))=r^q\lambda^d(B(1,0))$ implies that
$\lambda^d\bigl(\{x\in\rd :
\tau(x)=r\}\bigr)=0$, and it follows that $L_f$ is also continuous at point $1$ and thus on $(0,+\infty)$. Moreover, for any $r>0$ we have
$L_f(r^{-1})=-L_f(r)$. When $rs\geq 1$ with $r,s>0$ a change of
variables yields
$$
L_f(rs)=\int_{1\leq\tau(x)\leq
rs}f(x)\,dx=\int_{1\leq\tau(s^Ey)\leq rs}f(s^Ey)s^q\,dy
=\int_{s^{-1}\leq\tau(y)\leq r}f(y)\,dy.$$
Let us assume for instance that $1\le s^{-1}\le r$. Then, by continuity of $L_f$,
$$\int_{s^{-1}\leq\tau(y)\leq r}f(y)\,dy=\int_{1\leq\tau(y)\leq
r}f(y)\,dy-\int_{1\leq\tau(y)\le s^{-1}}f(y)\,dy=L_f(r)-L_f(s^{-1}).$$
It follows using $L_f(s^{-1})=-L_f(s)$ that
\begin{equation}\label{additivite}
L_f(rs)=L_f(r)+L_f(s).
\end{equation}
Similarly we show that (\ref{additivite}) holds 
 for $s^{-1}\le 1\le r$ and $ s^{-1}\le r\le 1$ and thus for all $rs\ge 1$. Using again the fact that $L_f(r^{-1})=-L_f(r)$, for all $r>0$, (\ref{additivite}) is valid for
 all $r,s>0$.
By continuity of $L_f$ it follows that $L_f(r)=L_f(e)\log r$.
We set $\mu_f=L_f(e)$. If $g(r)=1_{]a,b]}(r)$ for some $0<a<b$ we
get
\begin{eqnarray*}
\int_\rd f(x)g(\tau(x))\,dx&=&\int_{a<\tau(x)\leq
b}f(x)\,dx=L_f(b)-L_f(a)\\&=&\mu_f\bigl(\log b-\log
a\bigr)=\mu_f\int_0^\infty g(r)\,\frac{dr}r .\end{eqnarray*} The
general result follows by taking linear combinations and limits of
these functions in the standard way.
\end{proof}

\begin{proof}[Proof of Proposition \ref{prop12}]
When $f\in C(S_0)$ define $\tilde f$ on $\Gamma$ by $\tilde
f(x)=\tau(x)^{-q}f\bigl(l(x)\bigr)$. The function $\tilde f$
satisfies the hypothesis of Lemma \ref{lem13}. If $f\geq 0$ then
$\mu_{\tilde f}=L_{\tilde f}(e)=\int_{1\leq\tau(x)\leq
e}\tau(x)^{-q}f\bigl(l(x)\bigr)\,dx\geq 0$. Moreover $\mu_{a\tilde
f}=a\mu_{\tilde f}$, $\mu_{\tilde f+\tilde g}=\mu_{\tilde
f}+\mu_{\tilde g}$ and the mapping $f\mapsto\mu_{\tilde f}$ is
continuous. Hence this mapping is a positive linear functional on
$C(S_0)$. Therefore there exists a Radon measure $\sigma$ on $S_0$
such that $\mu_{\tilde f}=\int_{S_0}f(\theta)\,\sigma(d\theta)$.

If $g_1\in C_c((0,\infty))$ we get from applying Lemma \ref{lem13} with $\tilde f$ and $g(r)=r^qg_1(r)$ that
\begin{equation*}
\begin{split}
\int_\rd f(l(x))g_1(\tau(x))\,dx &= \int_\rd\tilde
f(x)\tau(x)^qg_1(\tau(x))\,dx \\ &=\mu_f\int_0^\infty
g_1(r)r^{q-1}\,dr \\
&=\int_0^\infty\int_{S_0}f(\theta)\,\sigma(d\theta)\,g_1(r)r^{q-1}\,dr
.
\end{split}
\end{equation*}
Since linear combinations of functions of the form
$f(l(x))g_1(\tau(x))$ are dense in $L^1(\rd,dx)$ the result follows.
\end{proof}

\begin{cor}\label{cor14}
Let $\beta\in\rr$ and suppose $f : \rd\to\mathbb C$ is measurable
such that $|f(x)|=O(\tau(x)^\beta)$. If $\beta>-q$ then $f$ is
integrable near $0$, and if $\beta<-q$ then $f$ is integrable near
infinity.
\end{cor}

We are now in position to define the class of $E$-homogeneous
functions and a important subclass needed in the moving average
representation of OSSRFs. Let $E$ be a $d\times d$ matrix as above
such that $0<a_1<\dots<a_p$ and for $x\in\Gamma$ let
$\bigl(\tau(x),l(x)\bigr)$ be the {\it polar coordinates}
associated to $E$, that is $x=\tau(x)^El(x)$.

\begin{defi}\label{def15}
Let $\varphi : \rd\to\mathbb C$ be any function. We say that
$\varphi$ is $E$-{\it homogeneous} if $\varphi(c^Ex)=c\varphi(x)$
for all $c>0$ and $x\in\Gamma$.
\end{defi}

It follows that an $E$-homogeneous function $\varphi$ is
completely determined by its values on $S_0$, since
$\varphi(x)=\varphi(\tau(x)^El(x))=\tau(x)\varphi(l(x))$. Observe
that if $\varphi$  is $E$-homogeneous and
continuous with positive values on $\Gamma$, then
\begin{equation}\label{F2}
M_\varphi=\max_{\theta\in S_0}\varphi(\theta)>0\ \text{and}\
m_\varphi=\min_{\theta\in S_0}\varphi(\theta)>0 .
\end{equation}
Moreover by continuity we necessarily have $\varphi(0)=0$.

\begin{defi}\label{def16}
Let $\beta>0$. A continuous function $\varphi : \rd \to
[0,\infty)$ is called {\it $\left(\beta, E\right)$-admissible}, if
$\varphi(x)>0$ for all $x\neq 0$ and for any $0<A<B$ there exists
a positive constant $C>0$ such that, for $A\le \|y\|\le B$,
 \[ \tau(x)\le 1\Rightarrow|\varphi(x+y)-\varphi(y)|\leq C\tau(x)^{\beta}.\]
\end{defi}

\begin{remark}\label{rem17}
If a continuous function $\varphi : \rd\to [0,\infty)$ is positive and Lipschitz  on $\Gamma$, that is $|\varphi(x)-\varphi(y)|\leq C\|x-y\|_0$ for
$x, y \in \Gamma$, then $\varphi$ is
$\left(\beta, E\right)$-admissible for all $\beta<a_1$ by Lemma \ref{lem10}.
\end{remark}

\begin{remark}
If $\varphi$ is $\left(\beta, E\right)$-admissible then $\beta\le
a_1$. In fact, if $\rd=V_1\oplus\cdots\oplus V_p$ is the spectral
decomposition of $\rd$ with respect to $E$ (see \cite{thebook},
Chapter 2 for details), by restricting the argument of the proof
of Lemma \ref{lem10} to the space $V_1$ on can show that for any
$\delta>0$ there exists a constant $C>0$ such that $\tau(x)\leq
C\|x\|_0^{1/a_1-\delta}$ for all $x\in V_1$ with $\|x\|_0\leq 1$.
Then, if for some fixed nonzero $u\in V_1$ we consider the
function $t\mapsto\varphi(tu)$ we get for $\delta_1=\beta\delta$
that $|\varphi(tu+su)-\varphi(su)|\leq C|t|^{\beta/a_1-\delta_1}$
for all small $t$ and $s$ bounded away from zero and infinity. If
one would have $\beta>a_1$, one could chose $\delta>0$ such that
$\beta/a_1-\delta_1>1$ and hence there would exist a constant
$K>0$ such that $\varphi(tu)=K$ for all $t\neq 0$. But since
$\varphi$ is continuous and $\varphi(0)=0$ this is impossible.
\end{remark}

\begin{remark}\label{expsym}
In general the exponent $E$ of a homogeneous function $\varphi$ is not unique.  It is easy to check that $\varphi(x)\to\infty$ as $\|x\|\to\infty$, and then Theorem 5.2.13 in \cite{thebook} implies that the set of possible exponents is $E+T\mathcal{S}(\varphi)$ where $E$ is any exponent, $\mathcal{S}(\varphi)$ is the set of symmetries of $\varphi$, and $T\mathcal{S}(\varphi)$ is the tangent space at the identity.  Here we say that $A$ is a symmetry of $\varphi$ if $\varphi(Ax)=\varphi(x)$ for all $x\in\rd$.  The symmetries $\mathcal{S}(\varphi)$ form a Lie group, and the tangent space consists of all derivatives $x'(0)$ of smooth curves $x(t)$ on $\mathcal{S}(\varphi)$ for which $x(0)=I$ the identity.  For example, if $\varphi$ is rotationally invariant then $\mathcal{S}(\varphi)$ is the orthogonal group and $T\mathcal{S}(\varphi)$ is the linear space of skew-symmetric matrices.  Although exponents are not unique, Theorem 5.2.14 in \cite{thebook} shows that every exponent $E$ of a homogeneous function $\varphi$ has the same real spectrum $0<a_1<\cdots<a_p$ and induces the same spectral decomposition $\rd=V_1\oplus\ldots\oplus V_p$, since these structural components describe the growth properties of the homogeneous function.  In particular, the function $r\mapsto \varphi(rx)$ grows like $r^{1/a_i}$ for any nonzero $x\in V_i$, see Section 5.3 in \cite{thebook} for more details.
\end{remark}

We conclude this section by examples of $\left(\beta,
E\right)$-admissible, $E$-homogeneous functions $\varphi :\rd\to[0,\infty)$ used in Theorem \ref{thm21} below to define a
moving average representation of OSSRFs
$\{X_\varphi(x)\}_{x\in\rd}$. Let us denote $<.,.>$ the standard inner product on $\rd$
and $E^t$ the transpose of any $d\times d$-matrix $E$ with respect to this inner product. The following class of examples is
inspired by the log-characteristic function of a full operator
stable law on $\rd$. See \cite{thebook} for details.

\begin{theorem}\label{thm41}
Assume $E$ is a real $d\times d$-matrix such that the real parts of the eigenvalues satisfy $1/2<a_1<a_2\ldots <a_p$ for $p\le d$. Assume $M(d\theta)$ is a
finite measure on the unit sphere $S_0$ corresponding to $E$ such
that
\[ \spann\bigl\{r^{E^t}\theta : r>0, \theta\in\supp(M)\bigr\}=\rd .\]
Then
\[
\varphi(x)=\int_{S_0}\int_0^\infty\bigl(1-\cos\bigl(\ip{x}{r^{E^t}\theta}\bigr)\bigr)\,\frac{dr}{r^2}\,M(d\theta)\]
is a continuous,  $E$-homogeneous function such that
$\varphi(x)>0$ for all $x\in\Gamma$. Moreover
$\varphi$ is $\left(\beta, E\right)$-admissible for
 $\beta<\min{\left(a_1,\frac{a_1}{a_p}\right)}$ if $a_1\le 1$ and $\beta=1$ if $a_1>1$.
\end{theorem}

\begin{proof}
Let $a_1>1/2$ denote the smallest real part of the eigenvalues of
$E$. Since $E$ and $E^t$ have the same eigenvalues, it follows from
Theorem 2.2.4 of \cite{thebook}
 that  for any $\delta>0$ there exists a
constant $C>0$ such that $\|r^{E^t}\theta\|_0\leq Cr^{a_1-\delta}$
for all $0<r\leq 1$ and $\theta_0\in S_0$. Therefore,  from
dominated convergence,  $\varphi$ is well-defined and
continuous on $\rd$. Moreover we have $\varphi(x)\geq 0$ and
$\varphi(x)=0$ implies $x=0$. A simple change of variable shows
that $\varphi(c^Ex)=c\varphi(x)$ for all $c>0$ and $x\in\rd$. It
remains to show that $\varphi$ is $(\beta,E)$-admissible. Using
the trigonometric identity
$\cos(a)-\cos(b)=-2\sin((a+b)/2)\sin((a-b)/2)$ we have for any
$x,y\in\rd$ that
\begin{equation}\label{accroissements}
|\varphi(x+y)-\varphi(y)|\leq
2\int_{S_0}\int_0^\infty\Bigl|\sin\Bigl(\frac{\ip{x+2y}{r^{E^t}\theta}}{2}\Bigr)
\sin\Bigl(\frac{\ip{x}{r^{E^t}\theta}}{2}\Bigr)\Bigr|\,\frac{dr}{r^2}\,M(d\theta).
\end{equation}
First, let us assume that $a_1>1$, then an upper bound of
\eqref{accroissements} is given by
$$
2\int_{S_0}\int_0^\infty\Bigl|\sin\Bigl(\frac{\ip{x}{r^{E^t}\theta}}{2}\Bigr)\Bigr|\,\frac{dr}{r^2}\,M(d\theta),$$
which is finite  because $a_1>1$, using $\|r^{E^t} \theta\|_0\leq Cr^{a_1-\delta}$ for all $0<r\leq 1$ and $\theta_0 \in S_0$, and elementary estimates.
Moreover writing $x=\tau(x)^El(x)$ a change of variables yields to
$$2\int_{S_0}\int_0^\infty\Bigl|\sin\Bigl(\frac{\ip{x}{r^{E^t}\theta}}{2}\Bigr)\Bigr|\,\frac{dr}{r^2}\,M(d\theta)=2\tau(x)\int_{S_0}\int_0^\infty\Bigl|\sin\Bigl(\frac{\ip{l(x)}{r^{E^t}\theta}}{2}\Bigr)\Bigr|\,\frac{dr}{r^2}\,M(d\theta),$$
which proves that $\varphi$ is 1-admissible.\\ Let us now consider
the case where $a_1\le 1$. Choose $\delta>0$ small enough. On one hand, for  $r\le
1$, one can find $C>0$ such that
$$\Bigl|\sin\Bigl(\frac{\ip{x+2y}{r^{E^t}\theta}}{2}\Bigr)
\sin\Bigl(\frac{\ip{x}{r^{E^t}\theta}}{2}\Bigr)\Bigr|\le
C\Bigl(2\| y\|_0 +\| x\|_0\Bigr)\| x\|_0r^{2a_1-\delta}.$$
On the other hand, it follows from  Theorem 2.2.4 of \cite{thebook}
that one can find $C>0$ such that
$\|r^{E^t} \theta\|_0\leq Cr^{a_p+\delta}$ for all $r\geq 1$ and $\theta_0 \in S_0$.
Thus, for
$\gamma<\min{\left(1,\frac{1}{a_p}\right)}$, using
$|\sin(u)|\leq |u|^\gamma$,  one can find $C>0$ such that
$$\Bigl|\sin\Bigl(\frac{\ip{x+2y}{r^{E^t}\theta}}{2}\Bigr)
\sin\Bigl(\frac{\ip{x}{r^{E^t}\theta}}{2}\Bigr)\Bigr|\le C\|
x\|_0^{\gamma}r^{\gamma a_p+\gamma\delta}.$$ Therefore,
by  substituting these  upper bounds into the right-hand side of (\ref{accroissements}) and integrating, for  some constant
$C>0$  we have shown that $|\varphi(x+y)-\varphi(y)|\leq
C\|x\|_0^\gamma$ for all $\|x\|_0\leq 1$ and $A\leq\|y\|_0\leq B$.

Since by Lemma \ref{lem10} $\|x\|_0\leq C\tau(x)^{a_1-\delta}$ for
$\tau(x)\leq 1$, the assertion follows with $\beta=\gamma(a_1- \delta)$.
\end{proof}

The following result gives a constructive description of a large
class of continuous, admissible $E$-homogeneous functions.

\begin{cor}\label{cor42}
Let $\theta_1,\dots,\theta_d$ be any basis of $\rd$, let
$0<\lambda_1\le\dots\le\lambda_d$ and $C_1,\dots,C_d>0$. Choose a
$d\times d$ matrix $E$ such that $E^t\theta_j=\lambda_j\theta_j$
for $j=1,\dots,d$. Then for any $\rho>0$, if $\rho<{2}\lambda_1$
the function
\[\varphi(x)=\Bigl(\sum_{j=1}^dC_j|\ip{x}{\theta_j}|^{\rho/\lambda_j}\Bigr)^{1/\rho}\]
is a continuous $E$-homogeneous and $\left(\beta,
E\right)$-admissible function for
$\beta<\min{\left({\lambda_1},{\rho}\frac{\lambda_1}{\lambda_d}\right)}$ if $\lambda_1\le\rho$ and $\beta={\rho}$ if $\lambda_1>\rho$.\\
\end{cor}

\begin{proof}
First observe that since $r^{E^t}\theta_j=r^{\lambda_j}\theta_j$
it follows that $\varphi(c^Ex)=c\varphi(x)$. Moreover $\varphi$ is
continuous. Let $B>A>0$, since $y\mapsto \sum_{j=1}^d
C_j|\ip{y}{\theta_j}|^{\rho/\lambda_j}$ is continuous and positive on $\Gamma$, by the mean value theorem, for $A\le\|y\|\le B$ and $\|x\|\le A/2$,
one can find $C>0$ such that
\begin{equation}\label{FE1}
|\varphi(x+y)-\varphi(y)|\leq
C\Bigl|\sum_{j=1}^dC_j|\ip{x+y}{\theta_j}|^{\rho/\lambda_j}-\sum_{j=1}^d
C_j|\ip{y}{\theta_j}|^{\rho/\lambda_j}\Bigr|.
\end{equation}
 Hence it remains to show that the
right hand side of \eqref{FE1} is $\left(\beta, E\right)$-admissible.
Let $M=\sum_{j=1}^d\gamma_j\varepsilon_{\theta_j}$
for suitable $\gamma_j>0$, where $\varepsilon_\theta$ denotes the
dirac mass in $\theta$. Let us define for $x\in \rd$,
$$\psi(x)=\int_{S_0}\int_0^\infty\bigl(1-\cos\bigl(\ip{x}{r^{(1/\rho)E^t}\theta}\bigr)\bigr)\,\frac{dr}{r^2}\,M(d\theta),$$
which is well defined
since $\rho<{2}\lambda_1$. Moreover, by Theorem \ref{thm41}, $\psi$ is $\left(\beta, (1/\rho)E\right)$-admissible for
$\beta<\min{\left(\frac{\lambda_1}{\rho},\frac{\lambda_1}{\lambda_d}\right)}$
if $\lambda_1<\rho$ and $\beta=1$ if $\lambda_1>\rho$. Let
$\tau_\rho(x)$ denote the radial part with respect to $(1/\rho)E$.
Then uniqueness implies that the radial part $\tau(x)$ with
respect to $E$ is given by $\tau(x)=\tau_\rho(x)^{1/\rho}$. Hence
$\psi$ is
 $\left(\beta, E\right)$-admissible for
 $\beta<\min{\left({\lambda_1},{\rho}\frac{\lambda_1}{\lambda_d}\right)}$ if $\lambda_1<\rho$ and $\beta={\rho}$ if $\lambda_1>\rho$.\\
Moreover, since
$r^{(1/\rho)E^t}\theta_j=r^{\lambda_j/\rho}\theta_j$ we get
\begin{equation*}
\begin{split}
\psi(x)=&\sum_{j=1}^d\gamma_j\int_0^\infty\bigl(1-\cos\bigl(r^{\lambda_j/\rho}|\ip{x}{\theta_j}|\bigr)\bigr)\,\frac{dr}{r^2}\\
=&\sum_{j=1}^d\frac{\rho\gamma_j}{\lambda_j}\Bigl(\int_0^\infty\bigl(1-\cos(s)\bigr)s^{-(\rho/\lambda_j)-1}\,ds\Bigr)
|\ip{x}{\theta_j}|^{\rho/\lambda_j} \\
=&\sum_{j=1}^dC_j|\ip{x}{\theta_j}|^{\rho/\lambda_j}.
\end{split}
\end{equation*}
This completes the proof.
\end{proof}

\section{Moving average representation}

In this section we consider a moving average representation of
OSSRFs and derive its basic properties. We first give sufficient
conditions such that the integral representation exists.  More
precisely, for $0<\alpha\leq 2$ we consider $Z_\alpha(dy)$ an
independently scattered $S\alpha S$ random measure on $\rd$ with
Lebesgue control measure $\lambda^d$. Then we define a moving
average representation of OSSRFs  using the basic fact that a random
integral $\int_\rd f(y)\,Z_\alpha(dy)$ exists if and only if
$\int_\rd|f(y)|^\alpha\,dy<\infty$.

Throughout this section we fix a real $d\times d$ matrix $E$ with
$0<a_1<\dots<a_p$ denoting the real parts of the eigenvalues of
$E$. As before, let $q=\trace(E)$.

\begin{theorem}\label{thm21}
Let $\beta>0$. Let $\varphi : \rd\to[0,\infty)$ be an
$E$-homogeneous, $\left(\beta, E\right)$-admissible function. Then
for any $0<\alpha\leq 2$ and any $0<H<\beta$ the random field
\begin{equation}\label{MAfield}
X_\varphi(x)=\int_\rd\Bigl(\varphi(x-y)^{H-q/\alpha}-\varphi(-y)^{H-q/\alpha}\Bigr)\,Z_\alpha(dy)\ ,x\in\rd
\end{equation}
exists and is stochastically continuous.
\end{theorem}

\begin{proof}
Let us recall that $X_\varphi(x)$ exists if and only if
$$\Gamma_{\varphi}^{\alpha}(x)=\int_\rd\Bigl|\varphi(x-y)^{H-q/\alpha}-\varphi(-y)^{H-q/\alpha}\Bigr|^\alpha\,dy<\infty.$$
Let us assume that $H\in (0,\beta)$. Observe that by \eqref{F2} and the fact that $\varphi$ is 
$E$-homogeneous,
$\varphi(z)\leq M_\varphi\tau(z)$ and $\varphi(z)\geq
m_\varphi\tau(z)$ for all $z\neq 0$. Fix any $x\in\Gamma$. Then,
\[
\Bigl|\varphi(x-y)^{H-q/\alpha}-\varphi(-y)^{H-q/\alpha}\Bigr|^\alpha\leq
C\Bigl(\tau(x-y)^{\alpha H-q}+\tau(y)^{\alpha H-q}\Bigr) .\] 
But for any
$R>0$ it
follows from Corollary \ref{cor14} that $\int_{\tau(y)\leq
R}\tau(y)^{\alpha H-q}\,dy<\infty$ if $H>0$. Moreover, by Lemma
\ref{lem11} $\{y : \tau(x-y)\leq R\}\subset\{y : \tau(y)\leq
K(R+\tau(x))\}$ and hence, by a change of variable we obtain using
Corollary \ref{cor14} again that, if $H>0$
\[\int_{\tau(y)\leq R}\tau(x-y)^{\alpha
H-q}\,dy=\int_{\tau(x-y)\leq R}\tau(y)^{\alpha H-q}\,dy\leq
\int_{\tau(y)\leq K(R+\tau(x))}\tau(y)^{\alpha H-q}\,dy<\infty. \]
\\ It remains to show that for some $R=R(x)>0$ we have
\begin{equation}\label{tt1}
\int_{\tau(y)>
R}\Bigl|\varphi(x+y)^{H-q/\alpha}-\varphi(y)^{H-q/\alpha}\Bigr|^\alpha\,dy<\infty.
\end{equation}
Observe that for $\tau(y)> R$, $\varphi(y)>0$, so we can write
$$\varphi(x+y)=\varphi\left(\varphi(y)^E\left(\varphi(y)^{-E}x+\varphi(y)^{-E}y\right)\right))=\varphi(y)
\varphi\left(\varphi(y)^{-E}x+\varphi(y)^{-E}y\right),$$ since
$\varphi$ is $E$-homogeneous. Moreover
$\varphi\left(\varphi(y)^{-E}y\right)=1$ and since $\varphi$ is
$\left(\beta,E\right)$ admissible, one can find $C>0$ such that
$$\left|\varphi\left(\varphi(y)^{-E}x+\varphi(y)^{-E}y\right)-1\right|\le C\tau\left(\varphi(y)^{-E}x\right)^{\beta}=C\varphi(y)^{-\beta}\tau\left(x\right)^{\beta}.$$
Hence by the mean value theorem applied to the function
$t^{H-q/\alpha}$ near $t=1$, one can find $C_1>0$ such that
\begin{equation*}
\begin{split}
\Bigl|\varphi(x+y)^{H-q/\alpha}-\varphi(y)^{H-q/\alpha}\Bigr| &
=\varphi(y)^{H-q/\alpha}\left|
\varphi(\varphi(y)^{-E}x+\varphi(y)^{-E}y)^{H-q/\alpha}-1\right|\\
& \le C_1\varphi(y)^{H-\beta-q/\alpha}\tau\left(x\right)^{\beta},
\end{split}
\end{equation*}
for all $\tau(y)> R$, where $R>0$ is chosen sufficiently large so that
$C\varphi(y)^{-\beta}\tau\left(x\right)^{\beta}< 1/2$ for all $\tau(y)> R$.
But $\varphi(y)^{H-\beta-q/\alpha}\le C_2\tau(y)^{H-\beta-q/\alpha}$
and by Corollary \ref{cor14} we know that $\int_{\tau(y)\geq
R}\tau(y)^{\alpha H-q-\alpha\beta}\,dy<\infty$ if $H<\beta$. This
allows  to conclude that  $\Gamma^{\alpha}_{\varphi}(x)$ is finite
for all $x\in \rd$. Let us now show that $X_\varphi$ is
stochastically continuous. Since $X_{\varphi}$ is a $S\alpha S$
field, it follows from Proposition 3.5.1 in \cite{sam:taqqu} that $X_{\varphi}$ is stochastically continuous if and only if, for all
$x_0\in \rd$,
$$\int_\rd\Bigl|\varphi(x_0+x-y)^{H-q/\alpha}-\varphi(x_0-y)^{H-q/\alpha}\Bigr|^\alpha\,dy\rightarrow
0\quad\text{as $x\to 0$.}$$ By a change a variables, this holds if
and only if
\begin{equation}\label{scalphama}
\Gamma_\varphi^{\alpha}(x)\rightarrow 0\quad\text{as $x\to 0$.}
\end{equation}
But $\varphi$ is continuous on $\rd$ so
$$\Bigl|\varphi(x-y)^{H-q/\alpha}-\varphi(-y)^{H-q/\alpha}\Bigr|^\alpha\rightarrow 0\quad\text{as $x\to 0$}$$ for almost every
 $y\in\rd$. Moreover, arguing as above, as soon as $\tau(x)\le 1$,
 for suitable $R>0$,
 one can find $C>0$ such that
\begin{multline*}
\Bigl|\varphi(x-y)^{H-q/\alpha}-\varphi(-y)^{H-q/\alpha}\Bigr|^\alpha
\\ \le
 C\left(\tau(y)^{\alpha H-q}\mathbf{1}_{\tau(y)\le K(R+1)}(y)+\tau(y)^{\alpha (H-\beta)-q}
 \mathbf{1}_{\tau(y)\ge R}(y)\right),
\end{multline*}
where $\mathbf 1_B(y)$ denotes the indicator function of a set
$B$.
Then (\ref{scalphama}) holds using dominated convergence.

 This concludes the proof.
\end{proof}

\begin{cor}\label{cor22}
Under the conditions of Theorem \ref{thm21}, the random field
$\{X_\varphi(x)\}_{x\in\rd}$ has the following properties:
\begin{enumerate}
\item[(a)] {\it operator scaling}, that is, for any $c>0$,
\begin{equation}\label{F4}
\{X_\varphi(c^Ex)\}_{x\in\rd}\eqfd\{c^{H}X_\varphi(x)\}_{x\in\rd}.
\end{equation}
\item[(b)] {\it stationary increments}, that is, for any $h\in\rd$,
\begin{equation}\label{F5}
\{X_\varphi(x+h)-X_\varphi(h)\}_{x\in\rd}\eqfd\{X_\varphi(x)\}_{x\in\rd}.
\end{equation}
\end{enumerate}
\end{cor}

\begin{proof}
We will only prove part (a). The proof of part (b) is left to the
reader. Fix any $x_1,\dots,x_m\in\rd$. Then \eqref{F4} follows if we
can show that for any $t_1,\dots,t_m\in\rr$ we have
\[\sum_{j=1}^mt_jX_\varphi(c^Ex_j)\eqd
c^{H}\sum_{j=1}^mt_jX_\varphi(x_j) .\] By a change of variable
together with $\varphi(c^Ex)=c\varphi(x)$ and the fact that $Z_\alpha(c^E dz)\eqd c^{q/\alpha}Z_\alpha(dz)$ we get
\begin{equation*}
\begin{split}
\sum_{j=1}^mt_jX_\varphi(c^Ex_j)  &=
\int_\rd\sum_{j=1}^mt_j\Bigl(\varphi(c^Ex_j-y)^{H-q/\alpha}
-\varphi(-y)^{H-q/\alpha}\Bigr)\,Z_\alpha(dy) \\ &\eqd
c^{q/\alpha}\int_\rd\sum_{j=1}^mt_j\Bigl(\varphi\bigl(c^E(x_j-z)\bigr)^{H-q/\alpha}
-\varphi(-c^Ez)^{H-q/\alpha}\Bigr)\,Z_\alpha(dz) \\
&=c^{H}\sum_{j=1}^mt_jX_\varphi(x_j)
\end{split}
\end{equation*}
and the proof is complete.
\end{proof}

\begin{remark}\label{iso}
Theorem \ref{thm21} and Corollary \ref{cor22} include the
following classical isotropic random fields as special cases.
Assume $\varphi(x)=\|x\|$ and $E=I$, the identity matrix. Observe
that $\varphi$ is an $E$-homogeneous, $(1,E)$-admissible function.
Then
\[
X_\varphi(x)=\int_\rd\bigl(\|x-y\|^{H-d/\alpha}-\|y\|^{H-d/\alpha}\bigr)\,Z_\alpha(dy)\]
Especially, if $\alpha=2$, then $\{X_\varphi(x)\}_{x\in\rd}$ is
known as the {\it L\'evy fractional Brownian field}. Note that in
this case, for any $0<\alpha\leq 2$ equation \eqref{F4} reduces to
the well-known self-similarity property
$\{X_\varphi(cx)\}_{x\in\rd}\eqfd\{c^HX_\varphi(x)\}_{x\in\rd}$.
Moreover our results  also include the well known one-dimensional
case $d=1$ of linear fractional stable motions and especially the
fractional Brownian motion when $\alpha=2$.
\end{remark}

\section{Harmonizable representation}
In this section we consider an harmonizable representation  of
OSSRFs and derive its basic properties. We first give necessary
and sufficient conditions such that the integral representation
exists and yields a stochastically continuous field. For
$0<\alpha\le 2$, let $W_{\alpha}(d\xi)$ be a complex isotropic
$S\alpha S$ random measure with Lebesgue control measure (see \cite{sam:taqqu} p.\ 281).

Throughout this section we fix a real $d\times d$ matrix $E$ with
$0<a_1<\dots<a_p$ denoting the real parts of the eigenvalues of
$E$. As before, let $q=\trace(E)$.

\begin{theorem}\label{thharm}
Let $\psi : \rd\to[0,\infty)$ be a continuous, $E^t$-homogeneous
function such that $\psi(x)\neq 0$ for $x\neq 0$. Then for any
$0<\alpha\leq 2$ the random field
\begin{equation}\label{Hfield}
X_\psi(x)=\Re\int_\rd\Bigl(e^{i<x,\xi>}-1\Bigr)\psi(\xi)^{-H-q/\alpha}\,W_\alpha(d\xi)\ ,x\in\rd
\end{equation}
exists and is stochastically continuous if and only if $H\in
(0,a_1)$.
\end{theorem}

\begin{proof}
Let us recall that $X_\psi(x)$ exists if and only if
$$\Gamma_\psi^{\alpha}(x):=\int_\rd\left|e^{i<x,\xi>}-1\right|^{\alpha}\psi(\xi)^{-\alpha H-q}\,d\xi <+\infty.$$
Let us assume that $H\in (0,a_1)$. By integration in {\it polar
coordinates} for $E^t$ given by Proposition \ref{prop12},
$$\Gamma_\psi^{\alpha}(x)=
\int_0^\infty\int_{S_0}\left|e^{i<x,r^{E^t}\theta>}-1\right|^{\alpha}r^{-\alpha
H-1}\psi(\theta)^{-\alpha H-q}\,\sigma(d\theta)\,\,dr.$$ For
$\delta \in (0, H-a_1)$, by considering the cases $r>1$ and $0\le r\le 1$ separately and using the same spectral bounds on the growth of $\| r^{E^t}\|$ as in the proof of Lemma \ref{lem10}, one
can find $C>0$ such that
$$\left|e^{i<x,r^{E^t}\theta>}-1\right|^{\alpha}\le C\left(1+\| x\|^{\alpha}\right)\min{(r^{\alpha (a_1-\delta)}, 1}).$$
Moreover, since $\psi$ is continuous with positive values on the sphere $S_0$, and hence bounded
away from zero,
$$\int_{S_0}\psi(\theta)^{-\alpha H-q}\,\sigma(d\theta)<\infty.$$
This allows  us to conclude that  $\Gamma_\psi^{\alpha}(x)$ is finite
for all $x\in \rd$. Let us show now that $X_\psi$ is
stochastically continuous. Since $X_{\psi}$ is a $S\alpha S$
field, it is stochastically continuous if and only if, for all
$x_0\in \rd$,
$$\int_\rd\Bigl|\left(\left(e^{i<x_0+x,\xi>}-1\right)-\left(e^{i<x_0,\xi>}-1\right)\right)
\Bigr|^\alpha\,\psi(\xi)^{-\alpha H-q}\,d\xi\to 0\quad\text{as
$x\to 0$}$$ that is, equivalently,
\begin{equation}\label{scalpha}
\Gamma_\psi^{\alpha}(x)\to 0\quad\text{as $x\to 0$.}
\end{equation}
It is straightforward to see that (\ref{scalpha}) holds for $H\in (0,a_1)$, using dominated convergence and the upper bound computed above.\\
Conversely, let us assume that $X_{\psi}$ exists and that it is
stochastically continuous. Let us remark that in this case
$\Gamma_\psi^{\alpha}(x)$ exists for all $x\in \rd$ and satisfies,
for all $\lambda>0$
$$\Gamma_\psi^{\alpha}(\lambda^Ex)=\lambda^{\alpha H}\Gamma_\psi^{\alpha}(x).$$
Let us fix any $x\in\rd$, with $x\neq 0$ and let us notice that
$\Gamma_\psi^{\alpha}(x)\neq 0$. Since $X_{\psi}$ is
stochastically continuous, by (\ref{scalpha}) 
$$\lambda^{\alpha H}\Gamma_\psi^{\alpha}(x)\rightarrow0\quad\text{as $\lambda\to 0$,}$$
which implies that $H>0$.\\
Let us now prove that $H<a_1$.\\
\underline{First case}: Assume that $a_1$ is an eigenvalue of $E$.
Then there exist $\theta_1\in \rd$ such that $\|\theta_1\|=1$ and
$E\theta_1=a_1\theta_1$. Therefore
$$\Gamma_\psi^{\alpha}(\theta_1)=\int_0^\infty\int_{S_0}\left|e^{i<\theta_1,r^{E^t}\theta>}-1\right|^{\alpha}r^{-\alpha H-1}\psi(\theta)^{-\alpha H-q}\,\sigma(d\theta)\,\,dr,$$
with
$$\bigl|<\theta_1,r^{E^t}\theta>\bigr|=r^{a_1}|<\theta_1,\theta>|\le
Cr^{a_1}.$$ Then, for $r\le \left(\frac{\pi}{C}\right)^{1/{a_1}}$,
$$\bigl|e^{i<\theta_1,r^{E^t}\theta>}-1\bigr|=2\bigl|\sin{\bigl(\frac{<\theta_1,r^{E^t}\theta>}{2}\bigr)}\bigr|\ge 2r^{a_1}\frac{|<\theta_1,\theta>|}{\pi},$$
and hence
$$\Gamma_\psi^{\alpha}(\theta_1)\ge\frac{1}{\pi^{\alpha}}\int_0^{\left(\frac{\pi}{C}\right)^{1/{a_1}}}\int_{S_0}\left|<\theta_1,\theta>\right|^{\alpha}r^{-\alpha (H-a_1)-1}\psi(\theta)^{-\alpha H-q}\,\sigma(d\theta)\,\,dr.$$
Since $\psi$ is positive on the sphere $S_0$,
$$\int_{S_0}\left|<\theta_1,\theta>\right|^{\alpha}\psi(\theta)^{-\alpha H-q}\,\sigma(d\theta)>0,$$
and then $\Gamma_\psi^{\alpha}(\theta_1)<+\infty$ implies that $H<a_1$.\\
\underline{Second case}: Assume that $a_1$ is not an eigenvalue of
$E$. Then there exists $b_1\in\rr$ such that $\lambda_1=a_1+ib_1$
and $\overline{\lambda_1}$ are complex eigenvalues of $E$. One can
find $\theta_1, \gamma_1\in\rd$, with
$\|\theta_1\|=\|\gamma_1\|=1$ such that
\begin{eqnarray*}
r^E\theta_1&=&r^{a_1}\left(\cos{\left(b_1\log r\right)}\theta_1+\sin\left(b_1\log r\right)\gamma_1\right)\\
r^E\gamma_1&=&r^{a_1}\left(-\sin\left(b_1\log r\right) \theta_1+\cos{\left(b_1\log r\right)}\gamma_1\right).
\end{eqnarray*}
Then it can be shown using the inequality $|e^{i\omega}-1|\geq |\omega|/\pi$ for $|\omega|<\pi$ that a lower bound of $\Gamma_\psi^{\alpha}(\theta_1)+\Gamma_\psi^{\alpha}(\gamma_1)$ is given by
$$
\frac{1}{\pi^{\alpha}}\int_0^{\left(\frac{\pi}{2C}\right)^{1/{a_1}}}\int_{S_0}\left(\left|<r^E\theta_1,\theta>\right|^{\alpha}+\left|<r^E\gamma_1,\theta>\right|^{\alpha}\right)r^{-\alpha
H-1}\psi(\theta)^{-\alpha H-q}\,\sigma(d\theta)\,\,dr.$$ Observe
that for $a,b\geq 0$ we have
$a^\alpha+b^\alpha\geq(a^2+b^2)^{\alpha/2}$. Therefore
\begin{eqnarray*}
\left|<r^E\theta_1,\theta>\right|^{\alpha}+\left|<r^E\gamma_1,\theta>\right|^{\alpha}&\ge&\left(\left|<r^E\theta_1,\theta>\right|^{2}+\left|<r^E\gamma_1,\theta>\right|^{2}\right)^{\alpha/2}\\
&\ge & r^{\alpha a_1}\left(\left|<\theta_1,\theta>\right|^{2}+\left|<\gamma_1,\theta>\right|^{2}\right)^{\alpha/2}.
\end{eqnarray*}
Then we conclude as in the first case that $H<a_1$. The proof is
complete.
\end{proof}

\begin{cor}\label{cor32}
Under the conditions of Theorem \ref{thharm}, the random field
$\{X_\psi(x)\}_{x\in\rd}$ has the following properties:
\begin{enumerate}
\item[(a)] {\it operator scaling}, that is, for any $c>0$,
\begin{equation}\label{Ff4}
\{X_\psi(c^Ex)\}_{x\in\rd}\eqfd\{c^{H}X_\psi(x)\}_{x\in\rd}.
\end{equation}
\item[(b)] {\it stationary increments}, that is, for any $h\in\rd$,
\begin{equation}\label{Ff5}
\{X_\psi(x+h)-X_\psi(h)\}_{x\in\rd}\eqfd\{X_\psi(x)\}_{x\in\rd}
\end{equation}
\end{enumerate}
\end{cor}

\begin{proof}
Let us recall that by corollary 6.3.2 of \cite{sam:taqqu}, for
$f\in L^{\alpha}(\rd)$, the characteristic function of the random
variable $Y=\Re\int_{\rd}f(y)W_\alpha(dy)$ is given by
\begin{equation}\label{characteristic function harm}
\Exp\left(e^{itY}\right)=\exp\Bigl(-c_0|t|^{\alpha}\int_{\rd}
|f(y)|^\alpha dy\Bigr)\quad\text{where}\quad c_0=\frac
1{2\pi}\int_0^\pi(\cos\theta)^2\,d\theta .
\end{equation}
Hence, for any
$x_1,\dots,x_m\in\rd$, the finite dimensional characteristic function
of $\left(X_{\psi}(x_1),\ldots,X_{\psi}(x_m)\right)$ is given by
$$\Exp\Bigl(\exp\Bigl(i\sum_{j=1}^mt_jX_{\psi}(x_j)\Bigr)\Bigr)=
\exp\Bigl(-c_0\int_{\rd}
\Bigl|\sum_{j=1}^mt_j\left(e^{i<x_j,\xi>}-1\right)\Bigr|^\alpha
\psi(\xi)^{-\alpha H-q}d\xi\Bigr),$$ for any
$t_1,\dots,t_m\in\rr$. Thus, for any $c>0$, by a change of
variable $\gamma=c^{E^t}\xi$ in the integral of the right side,
since $\psi$ is an $E^t$ homogeneous function, we get
$$\Exp\Bigl(\exp\Bigl(i\sum_{j=1}^mt_jX_{\psi}(c^Ex_j)\Bigr)\Bigr)=\Exp\Bigl(\exp\Bigl(i\sum_{j=1}^mt_jc^HX_{\psi}(x_j)\Bigr)\Bigr),$$
which proves (a). Furthermore, for any $h\in\rd$ and $x\in\rd$, we
have that
$$X_{\psi}(x+h)-X_{\psi}(x)=\Re\int_\rd e^{i<h,\xi>}\Bigl(e^{i<x,\xi>}-1\Bigr)\psi(\xi)^{-H-q/\alpha}\,W_\alpha(d\xi).$$ Hence
\begin{eqnarray*}
&&\Exp\Bigl(\exp\Bigl(i\sum_{j=1}^mt_j\left(X_{\psi}(x_j+h)-X_{\psi}(x_j)\right)\Bigr)\Bigr)\\&=&
\exp\Bigl(-c_0\int_{\rd}
\Bigl|\sum_{j=1}^mt_je^{i<h,\xi>}\left(e^{i<x_j,\xi>}-1\right)\Bigr|^\alpha
\psi(\xi)^{-\alpha H-q}d\xi\Bigr)\\
&=&\Exp\Bigl(\exp\Bigl(i\sum_{j=1}^mt_jX_{\psi}(x_j)\Bigr)\Bigr),
\end{eqnarray*}
proving (b).
\end{proof}

\begin{remark} In the Gaussian case, the covariance function of the random field $X_\varphi(x)$ defined by the moving average representation \eqref{MAfield} can be computed by an argument similar to Proposition 8.1.4 of \cite{sam:taqqu}.  Let $\sigma^2_\theta=\Exp[(X_\varphi(\theta))^2]$ for any unit vector $\theta$, and define $\tau(x)$ and $\ell(x)$ as before so that $x=\tau(x)^E\ell(x)$.  Using Corollary \ref{cor22} (a) it follows that $\Exp[(X_\varphi(x))^2]=\tau(x)^{2H}\sigma^2_{\ell(x)}$, and then we can use the fact that 
$2X_\varphi(x)X_\varphi(y)=X_\varphi(x)^2+X_\varphi(y)^2-(X_\varphi(x)-X_\varphi(y))^2$
to conclude that
\begin{equation}\label{ACF-MA}
\Exp\left[X_\varphi(x)X_\varphi(y)\right]=\tfrac 12\left[ \tau(x)^{2H}\sigma^2_{\ell(x)} +\tau(y)^{2H}\sigma^2_{\ell(y)}-\tau(x-y)^{2H}\sigma^2_{\ell(x-y)}\right].
\end{equation}
In the isotropic case discussed in Remark \ref{iso} we have $\tau(x)=\|x\|$ and $\ell(x)=x/\|x\|$, and a change of variables in \eqref{MAfield} shows that $\sigma^2_\theta\equiv\sigma^2$ is the same for any unit vector, using the fact that $\varphi(Rx)=\varphi(x)$ for any orthogonal linear transformation $R$ in this case.  Then \eqref{ACF-MA} reduces to the familiar autocovariance function for a fractional Gaussian random field.  A similar argument shows that the autocovariance function of the random field defined by the harmonizable representation \eqref{Hfield} is given by 
\begin{equation}\label{ACF-H}
\Exp\left[X_\psi(x)X_\psi(y)\right]=\tfrac 12\left[ \tau(x)^{2H}\omega^2_{\ell(x)} +\tau(y)^{2H}\omega^2_{\ell(y)}-\tau(x-y)^{2H}\omega^2_{\ell(x-y)}\right]
\end{equation}
where $\omega^2_\theta=\Exp[(X_\psi(\theta))^2]$.  For the isotropic case, where \eqref{Hfield} reduces to the harmonizable representation \eqref{eq13} for a fractional Gaussian field, we again note that $\omega^2_\theta$ is constant over the unit sphere.  Since a mean zero Gaussian random field is determined by its autocovariance function, we recover the well-known fact that the moving average and harmonizable representations of the fractional Gaussian random field differ by at most a constant factor.  It does not seem possible to extend this argument to the general case of operator scaling Gaussian random fields, since it would be difficult to compare $\sigma^2_\theta$ to $\omega^2_\theta$ in this case.  Hence there remains an interesting open question under which relationship
between the functions $\varphi$ and $\psi$ in the Gaussian case
the moving average representation of Theorem \ref{thm21} and the
harmonizable representation of Theorem \ref{thharm} are
equivalent.  
\end{remark}

\begin{remark}
Many random fields occurring in applications have Hurst indices that vary with coordinate \cite{benson,Bonami}.  Consider a random field satisfying \eqref{eq11}, and suppose that the matrix $E$ has an eigenvector $e$ with associated real eigenvalue $\lambda$.  Then it follows from \eqref{eq11} that the stochastic process $r\mapsto X(re)$ is self-similar with 
\[\{X(c^\lambda re)\}_{r\in\rr}\eqfd\{c^H X(re)\}_{r\in\rr}\quad\text{for all
$c>0$,}\]
so that the Hurst index of this process is $H/\lambda$.  If $E$ has a basis of eigenvectors with distinct real eigenvalues, then the projections of this random field onto the eigenvector directions yield processes with different Hurst indices in each coordinate.  This also shows that the usual methods for estimating the Hurst index, such as  rescaled range analysis \cite{hurst} and dispersional analysis  \cite{caccia}, can also be applied to estimate the scaling indices of the operator scaling random field from data, once the proper coordinates are established.  Estimating these coordinate directions from data is an interesting open question.  In some practical applications, these coordinates are known from the problem setup.  For example, in a groundwater aquifer the coordinates of the hydraulic conductivity field are thought to correspond to the vertical, the direction of horizontal mean flow, and the horizontal direction perpendicular to the mean flow \cite{benson}.  In fractured rock, the scaling coordinates of the transmissivity field correspond to the main fracture orientations, and are usually not mutually perpendicular \cite{Schumer}.  Similarly, in material science, the crack fronts determine the natural coordinates \cite{Ponson}.  We caution, however, that estimating the Hurst index in the wrong (non-eigenvalue) coordinates is likely to be misleading, because in those directions the field is not self-similar.  Finally, we note that the parameters $E,H$ in \eqref{eq11} are not unique.  If \eqref{eq11} holds, then we also have $\{X(c^{E'}x)\}\eqfd\{cX(x)\}$ where $E'=(1/H)E$, so that the Hurst indices of the random field are the ratio of $H$ and the eigenvalues of $E$, as already noted.  Furthermore, the exponents of an admissible function are not unique, because of possible symmetries, as discussed previously in Remark \ref{expsym}.  Hence the Hurst index of each component is really an estimate of $H/a_i$ where $0<a_1<\cdots<a_p$ is the real spectrum of $E$, and these indices, as well as the coordinate system in which they pertain, are the same for any choice of $H$ and $E$.  
\end{remark}

We have already seen that the OSSRFs, defined by a moving average
or a harmonizable representation were stochastically continuous.
In the next section we show that in the Gaussian case $\alpha=2$
one can get H{\"o}lder regularity for the sample paths.

\section{Gaussian OSSRFs}
In this section, we are interested  in the smoothness of the
sample paths of Gaussian OSSRFs given by Theorem \ref{thm21} or
Theorem \ref{thharm}, respectively. Moreover we compute the box-
and the Hausdorff-dimension of the graph of OSSRFs in these cases.
We follow the terminology used in \cite{Bonami}. Using their
definition of the H{\"o}lder
 critical exponent of a random process (Definition 5) we state the
 following definition.
 \begin{defi} Let $\gamma\in(0,1)$. A random field
 $\{X(x)\}_{x\in\rd}$ is said to have H\"older
 critical exponent $\gamma$ whenever it satisfies the following
two properties:
\begin{enumerate}
\item[(a)] for any $s\in (0,\gamma)$, the sample paths of $X$
satisfy almost surely a
 uniform H{\"o}lder condition of order $s$ on any compact set,
 that is
 for any compact set $K\subset\rd$, there exists a positive random
 variable $A$ such that
 $$\left|X(x)-X(y)\right|\le A\|x-y\|^{s}\quad\text{for all $x,y \in
 K$.}$$
\item[(b)] for any $s\in (\gamma,1)$, almost surely the sample paths of
$X$ fail to
 satisfy any uniform H{\"o}lder condition of order $s$.
\end{enumerate}
\end{defi}

For a Gaussian random field $X$  a well known result links the
H\"older regularity of the sample paths $x\mapsto X(x,\omega)$ to
those of the quadratic mean. Let us recall this property when the
field also has stationary increments. We refer to \cite{Adler}
Theorem 8.3.2 and Theorem 3.3.2 for a detailed proof.

\begin{prop}\label{regularity}
Let $\{X(x)\}_{x\in\rd}$ be a Gaussian random field with stationary
increments. Let $\gamma\in(0,1)$ and assume that
$$\gamma=\sup\left\{s>0;
\Exp\left(\left(X(x)-X(0)\right)^2\right)=o_{\| x\|\rightarrow 0}\left(\| x\|^{2s}\right)\right\}.$$
Then, for any $s\in (0,\gamma)$, any continuous version of $X$  satisfies almost surely a
 uniform H{\"o}lder condition of order $s$ on any compact set.\\
 If moreover
 $$\gamma=\inf\left\{s>0; \| x\|^{2s}
=o_{\| x\|\rightarrow
0}\left(\Exp\left(\left(X(x)-X(0)\right)^2\right)\right)\right\},$$
then any continuous version of $X$ admits $\gamma$ as the H\"older
 critical exponent.
\end{prop}

The previous definition and proposition are given in \cite{Bonami}
for random processes ($d=1$) in order to study regularity
properties of a field along straight lines. More precisely, when
$\{X(x)\}_{x\in\rd}$ is a random field, it is also interesting to
study the H\"older regularity of the process
$\{X(x_0+tu)\}_{t\in\rr}$, for $x_0\in \rd$ and $u$ a unit vector.
This will provide some additional directional regularity
information. For $\{X(x)\}_{x\in\rd}$ with stationary increments,
one only has to consider $\{X(tu)\}_{t\in\rr}$ for all directions
$u$. Let us recall Definition 6 of \cite{Bonami}.
\begin{defi} Let $\{X(x)\}_{x\in\rd}$ with stationary increments
and let $u$ be any direction of the unit sphere. If the process
$\{X(tu)\}_{t\in\rr}$ has H\"older critical exponent $\gamma(u)$ we
say that $X$ admits $\gamma(u)$ as directional regularity in
direction $u$.
 \end{defi}

Let us investigate these properties for the Gaussian OSSRFs given
by Theorem \ref{thm21} or Theorem \ref{thharm}, respectively.
Throughout this section we fix a real $d\times d$ matrix $E$ with
$0<a_1<\dots<a_p$ denoting the real parts of the eigenvalues of
$E$. Following \cite{thebook}, Section 2.1, let $V_1,\ldots,V_p$
be the spectral decomposition of $\rd$ with respect to $E$. For
$i=1,\ldots,p$, let us denote
$$W_i=V_1\oplus\ldots\oplus V_i,$$ and $W_0=\{0\}$. Observe that
$E|_{W_i}$ has $a_1<\dots<a_i$ as real parts of the eigenvalues.
As before let $q=\trace(E)$.

\begin{theorem}\label{holder}
Let $\varphi : \rd\to[0,\infty)$ be an $E$-homogeneous,
$(\beta,E)$-admissible function. For $0<H<\beta$ let
$X_\varphi$ be the moving average Gaussian OSSRF given by Theorem
\ref{thm21}. Moreover let $\psi : \rd\to[0,\infty)$ be a
continuous $E^t$-homogeneous function with $\psi(x)>0$ for all $x\neq 0$. For $0<H<a_1$ let
$X_\psi$ be the harmonizable Gaussian OSSRF given by Theorem
\ref{thharm}. Then any continuous version of $X_{\varphi}$ and
$X_{\psi}$, respectively, admits $H/a_p$ as H\"older critical
exponent. Moreover, for any $i=1,\ldots,p$, for any direction
$u\in W_i\setminus W_{i-1}$, the fields $X_{\varphi}$ and
$X_{\psi}$ admit $H/a_i$ as directional regularity in direction
$u$.
\end{theorem}

\begin{proof}
Let $x\in \rr^d$. With a little abuse of notation we write
$X_{\varphi/\psi}$ to indicate that we either consider $X_\varphi$
or $X_\psi$. Observe that $X_{\varphi/\psi}(0)=0$ and in order to
apply Proposition \ref{regularity} we define
\[\Gamma^2_{\varphi/\psi}(x)=\mathbf{E}\left(X_{\varphi/\psi}(x)^2\right)=\begin{cases}
\int_\rd\left|\varphi(x-y)^{H-q/2}-\varphi(-y)^{H-q/2}\right|^2dy  \\
4\int_\rd\sin^2\left(\frac{<x,\xi>}{2}\right)\psi(\xi)^{-2H-q}d\xi
\end{cases}\]
Using polar coordinates with respect to $E$, it is straightforward to see that
\begin{equation}\label{variance}
\Gamma^2_{\varphi/\psi}(x)=\tau(x)^{2H}\Gamma_{\varphi/\psi}(l(x)),
\end{equation}
where for all $\theta\in S_0$,
\begin{equation}\label{encsphere}
0<m\le \Gamma^2_{\varphi/\psi}(\theta)\le M,
\end{equation}
since $\Gamma^2_{\varphi/\psi}$ is continuous and positive on the compact set $S_0$.

For any $i=1,\dots,p$ let us fix $u\in W_i\setminus W_{i-1}$.
Since the spaces $V_1,\dots,V_p$ are $E$-invariant and the real
parts of the eigenvalues of $E|_{W_i}$ are $a_1<\dots<a_i$ it
follows as in the proof of Lemma \ref{lem10} by considering the
space $W_i$ instead of $\rd$, that for any small $\delta>0$ there
exists a constant $C_2=C_2(u)>0$ such that $\tau(tu)\leq
C_2|t|^{1/a_i-\delta}$ for any $|t|\leq 1$. Furthermore, observe
that if we write $u=u_i+\bar u_{i-1}$ with $u_i\in V_i$ and $\bar
u_{i-1}\in W_{i-1}$ we have $u_i\neq 0$. Writing
$tu=\tau(tu)^El(tu)$ and $l(tu)=l_i(tu)+\bar l_{i-1}(tu)$ with
$l_i(tu)\in V_i$ and $\bar l_{i-1}(tu)\in W_{i-1}$, it follows from
the $E$-invariance of the spectral decomposition that
$tu_i=\tau(tu)^El_i(tu)$ with $l_i(tu)\neq 0$. Since we have
$E=E_1\oplus\cdots\oplus E_p$ where every real part of the
eigenvalues of $E_i$ equals $a_i$ we conclude
\begin{equation*}
|t|\|u_i\|=\|\tau(tu)^El_i(tu)\|=\|\tau(tu)^{E_i}l_i(tu)\|\leq
\|\tau(tu)^{E_i}\|\|l_i(tu)\|\leq C\tau(tu)^{a_i-\delta}
\end{equation*}
for any $|t|\leq 1$ using the fact that $\|l_i(tu)\|\leq C_3$ for
any $|t|\leq 1$ and some $C_3>0$. Hence there exists a constant
$C_1=C_1(u)>0$ such that $\tau(tu)\geq C_1|t|^{1/a_i+\delta}$ for
any $|t|\leq 1$. Therefore we have shown that for all directions
$u\in W_i\setminus W_{i-1}$ and any small $\delta>0$ there exist
constants $C_1,C_2>0$, such that
\begin{equation}\label{sbound}
C_1| t|^{1/a_i+\delta}\le \tau(tu)\le C_2|
t|^{1/a_i-\delta}\quad\text{for all $|t|\leq 1$.}
\end{equation}
In view of \eqref{variance}, \eqref{encsphere} and \eqref{sbound}
we therefore get that for any direction $u\in W_i\setminus
W_{i-1}$ and any $\delta>0$ there exist constants $C_1,C_2>0$ such
that $C_1|t|^{2H/a_i+\delta}\leq \Gamma^2_{\varphi/\psi}(tu)\leq
C_2|t|^{2H/a_i-\delta}$ for $|t|\leq 1$, which by Proposition
\ref{regularity} shows that $X_{\varphi/\psi}$ admits $H/a_i$ as
directional regularity in direction $u$.
\\
It follows from this that for any $s\in (H/a_p,1)$ almost surely the sample
paths of $X_{\varphi/\psi}$ fail to satisfy any uniform H{\"o}lder
condition of order $s$, since $H/a_p$ is the H\"older
critical exponent of $X_{\varphi/\psi}$ in any direction of
$W_p\setminus W_{p-1}$. Finally, in view of \eqref{variance},
\eqref{encsphere} and Lemma \ref{lem10} we know that for any
$\delta>0$ there exists a constant $C>0$ such that
$\Gamma^2_{\varphi/\psi}(x)\leq C\|x\|^{2H/a_p-\delta}$ for
$\|x\|\leq 1$ and hence by Proposition \ref{regularity} it follows
that any continuous version of $X_{\varphi/\psi}$ satisfies almost surely a
uniform H\"older condition of order $s<H/a_p$ on any compact
set. This concludes the proof.
\end{proof}

Having described the H\"older regularity of Gaussian OSSRFs, a
natural question that arises is to determine the  box- and the
Hausdorff-dimensions of their graphs on a compact set. We refer to
Falconer \cite{Falconer} for the definitions and properties of
box- and the Hausdorff-dimension. Let us fix a compact set
$K\subset\rd$. For a random field $X$ on $\rd$ we consider
$\mathcal G(X)(\omega)=\{(x,X(x)(\omega)); x\in K\}$ the graph of
a realization of this field over the compact $K$. We will denote
$\mbox{dim}_{\mathcal H}\mathcal G(X)$, resp $\mbox{dim}_{\mathcal
B}\mathcal G(X)$, the Hausdorff-dimension and the box-dimension of
$\mathcal G(X)$, respectively.

It is a well understood fact that directional regularity implies
information about the Hausdorff-dimension of the field in that
direction. See e.g. \cite{Adler}, Chapter 8. As an immediate
corollary to Theorem \ref{holder} we get:
\begin{cor}\label{Hdir}
Under the assumptions of Theorem \ref{holder} we have for all
$i=1,\dots,p$ and all directions $u\in W_i\setminus W_{i-1}$ that
\[ \dim_{\mathcal H}\{(t,X_{\varphi/\psi}(tu)) :
t\in[0,1]\}=2-H/a_i\quad\text{a.s.}\]
\end{cor}

\begin{proof}
The result is a direct consequence of Theorem \ref{holder} and the
Corollary on page 204 of \cite{Adler}, using the fact that
$t\mapsto X_{\varphi/\psi}(tu)$ is a $\beta=H/a_i$-index process
and that $2-\beta\leq 1/\beta$ for $0<\beta<1$.
\end{proof}

Our next result investigates the global box- and
Hausdorff-dimension of Gaussian OSSRFs.

\begin{theorem}\label{XHD} Under the assumptions of Theorem \ref{holder}, for
any continuous version of $X_{\varphi}$ and $X_{\psi}$, almost
surely
$$\dim_{\mathcal H}\mathcal G(X_{\varphi/\psi})=\dim_{\mathcal B}\mathcal G(X_{\varphi/\psi}) = d+1-H/a_p.$$
\end{theorem}

\begin{proof}
Let us choose a continuous version of $X_{\varphi/\psi}$. From
Theorem \ref{holder}, for any $s<H/a_p$, the sample paths of
$X_{\varphi/\psi}$ satisfy almost surely a uniform H\"older condition of order
$s$ on $K$. Thus by a $d$-dimensional version of Corollary 11.2 of \cite{Falconer}, we
have
$$\mbox{ dim}_{\mathcal H}\mathcal G(X_{\varphi/\psi})\leq\overline{\mbox{ dim}_{\mathcal B}}\mathcal G(X_{\varphi/\psi}) \le d+1-s,\text{ a.s.}$$
where $\overline{\mbox{ dim}_{\mathcal B}}$ denotes the upper
box-dimension. Therefore
$$\mbox{ dim}_{\mathcal H}\mathcal G(X_{\varphi/\psi})\leq\overline{\mbox{ dim}_{\mathcal B}}\mathcal G(X_{\varphi/\psi}) \le d+1-H/a_p,\text{ a.s.}$$
and it remains to show that  a.s. $\mbox{ dim}_{\mathcal H}\mathcal
G(X_{\varphi/\psi})\ge d+1-H/a_p$.  Since the lower box
dimension satisfies $\underline{\mbox{ dim}_{\mathcal B}}\mathcal
G(X_{\varphi/\psi})\ge
\mbox{ dim}_{\mathcal H}\mathcal G(X_{\varphi/\psi})$ the proof is then complete.\\
We follow the same kind of ideas developed in \cite{Benassi} and
\cite{AyacheR}. Let $s>1$. Following the same argument as in Theorem 16.2 of \cite{Falconer}, in view of the Frostman criterion (Theorem 4.13 (a) in \cite{Falconer}), if one proves that the
integral $I_s$
$$I_s=\int_{K\times K}\mathbf{E}\left[\left((X_{\varphi/\psi}(x)-X_{\varphi/\psi}(y))^2+\| x-y\|^2\right)^{-s/2}\right]\,dx\,dy,$$
is finite, then almost surely $\mbox{ dim}_{\mathcal H}\mathcal
G(X_{\varphi/\psi})\ge s$.

As before, let $V_1,\dots,V_p$ denote the spectral decomposition
of $\rd$ with respect to $E$ and let $W_i=V_1+\cdots+V_i$. We will
choose an inner-product $(\cdot,\cdot)$ on $\rd$ which makes these
spaces mutually orthogonal and use the norm $\|x\|=(x,x)^{1/2}$.
Since all norms on $\rd$ are equivalent, this entails no loss of
generality.

Since by assumption $s>1$, the function $(\xi^2+1)^{-s/2}$ is in $
L^1(\rr)$ and its Fourier transform, denoted by $f_s$, is not only in $L^{\infty}(\rr)$ but also in
$L^1(\rr)$. Then we can write, using Fourier-inversion
(fundamental lemma in \cite{Benassi})
$$(\xi^2+1)^{-s/2}=\frac{1}{2\pi}\int_{\rr}e^{i\xi t}f_s(t)dt.$$
It follows that
\begin{eqnarray*}
& & \Exp\left[\left((X_{\varphi/\psi}(x)-X_{\varphi/\psi}(y))^2+\|
x-y\|^2\right)^{-s/2}\right]\\ 
&= & \frac{1}{2\pi}\|
x-y\|^{-s}\int_{\rr}\Exp\left(e^{it\frac{X_{\varphi/\psi}(x)-X_{\varphi/\psi}(y)}{
\| x-y\|}}\right)f_s(t)dt\\
&= & \frac{1}{2\pi}\| x-y\|^{-s}\int_{\rr}
e^{-\frac{t^2}{2}\frac{\Exp\left(X_{\varphi/\psi}(x)-X_{\varphi/\psi}(y))^2\right)}{\|
x-y\|^2}}f_s(t)dt,
\end{eqnarray*}
since $X_{\varphi/\psi}$ is Gaussian. Then, as $f_s\in L^{\infty}(\rr)$, one can find $C>0$ such that
\begin{eqnarray*}
& &\Exp\left[\left((X_{\varphi/\psi}(x)-X_{\varphi/\psi}(y))^2+\|
x-y\|^2\right)^{-s/2}\right]\\
& \le & C\|x-y\|^{1-s}\left(\Exp\left[\left(X_{\varphi/\psi}(x)-X_{\varphi/\psi}(y)\right)^2\right]\right)^{-1/2}\\
&\le &  Cm^{-1}\| x-y\|^{1-s}\tau(x-y)^{-H},
\end{eqnarray*}
according to (\ref{variance}) and (\ref{encsphere}) and using the fact that $X_{\varphi/\psi}$ has stationary increments.

Let us choose $A>0$ such that $K\subset \{x\in \rd; \|x\|\le
A/2\}$. Then for some constant $C>0$
$$I_s\le C\int_{\|x\|\le A}\| x\|^{1-s}\tau(x)^{-H}dx,$$
as long as the integral in the right hand side is bounded.\\
If $p=1$, by Lemma \ref{lem10}, for $\delta>0$, one can find $C>0$ such that, for $\|x\|\le A$,
$$\tau(x)^{-H}\le C\| x\|^{-H/a_p-\delta},$$
and hence $I_s$ is finite as soon as $s<d+1-H/a_p-\delta$.\\
If $p\geq 2$ let us write $x=x_p+y$ for some $x_p\in V_p$ and
$y\in W_{p-1}$ and write $x=\tau(x)^El(x)$ with $l(x)\in S_0$.
Decompose $l(x)=l_p(x)+\theta$ with $l_p(x)\in V_p$ and $\theta\in
W_{p-1}$. By the direct sum decomposition we see that
$x_p=\tau(x)^El_p(x)$ and $y=\tau(x)^E\theta$. Moreover, since
$V_p$ and $W_{p-1}$ are orthogonal in the chosen inner product it follows that $\|x\|\leq A$
implies both $\|x_p\|\leq A$ and $\|y\|\leq A$ in the associated norm. In view of the
proof of Lemma \ref{lem10}, restricted to the spaces $V_p$ and
$W_{p-1}$, respectively, it follows that for any $\delta>0$ and
some constants $C_1,C_2>0$, if $\|x\|\leq A$ then
$$\|x_p\|\le C_1\tau(x)^{a_p-\delta}\,\mbox{ and } \|y\|\le C_2\tau(x)^{a_1-\delta}.$$
Then one can find $c>0$ such that $$\tau(x)^H\ge c\| x_p\|^{H/a_p+\delta} \mbox{ and } \tau(x)^H\ge c\| y\|^{H/a_1+\delta}$$
and thus
$$\tau(x)^H\ge c/2\left(\| x_p\|^{H/a_p+\delta}+\| y\|^{H/a_1+\delta}\right).$$
Hence, for any $\delta>0$
$$I_s\le C\int_{\|x_p\|\le A}\int_{\|y\|\le A}\left(\| x_p\|^2+\| y\|^2\right)^{1/2-s/2}\left(\| x_p\|^{H/a_p+\delta}+\| y\|^{H/a_1+\delta}\right)^{-1}\,dy\,dx_p.$$
Let $k=\dim V_p$ and observe that in the present case $1\leq k\leq
d-1$.  By using polar coordinates for both $V_p$ and $W_{p-1}$,
for some constant $C>0$ we have  $I_s\le CJ_s$ where
$$J_s= \int_0^A\int_0^A(u^2+v^2)^{1/2-s/2}\left(u^{H/a_p+\delta}+v^{H/a_1+\delta}\right)^{-1}u^{k-1}v^{d-1-k}\,du\,dv.$$
The change of variables $u=tv$ yields
\begin{eqnarray*}
J_s&=&\int_0^A\int_0^{A/v}v^{d-s-H/a_p-\delta}(t^2+1)^{1/2-s/2}\left(t^{H/a_p+\delta}+v^{H/a_1-H/a_p}\right)^{-1}t^{k-1}\,dt\,dv\\
&\le &\left(\int_0^A v^{d-s-H/a_p-\delta}
dv\right)\left(\int_0^{+\infty}(t^2+1)^{1/2-s/2}t^{-H/a_p-\delta+k-1}dt\right).
\end{eqnarray*}
Since $\frac{H}{a_p}<1\le k$, the second term is bounded as soon
as $s>k+1-H/a_p-\delta$, whereas the first one is finite whenever
$s<d+1-H/a_p-\delta$. Thus, for all $\delta>0$ small enough, it follows that almost surely
$\mbox{ dim}_{\mathcal H}\mathcal G(X_{\varphi/\psi})\ge
d+1-H/a_p-\delta$ and the proof is complete.
\end{proof}

\begin{remark}
As pointed out in the introduction, the fractional Brownian sheet
$\{B_H(x)\}_{x\in\rd}$ is operator scaling with
$E=\diag(\alpha_1,\dots,\alpha_d)$ where $\alpha_i=1/H_i$, in fact
$\{B_H(c^Ex)\}_{x\in\rd}\eqfd\{cB_H(x)\}_{x\in\rd}$, but does not
have stationary increments. By Theorem 1.1 of \cite{ayache1} we
know that $\dim_{\mathcal H}\mathcal G(B_H)=\dim_{\mathcal
B}\mathcal G(B_H)=d+1-\min(\alpha_1^{-1},\dots,\alpha_d^{-1})$.
Now let $X_{\varphi/\psi}$ be our Gaussian OSSRFs considered
above. It follows from operator scaling, that
\[\{X_{\varphi/\psi}(c^{\frac
1HE}x)\}_{x\in\rd}\eqfd\{cX_{\varphi/\psi}(x)\}_{x\in\rd}\] too,
and that $\bar\alpha_i=a_i/H$, $i=1,\dots,p$ are the real parts of
the eigenvalues of $(1/H)E$. Hence Theorem \ref{XHD} can be
reformulated as
\[\dim_{\mathcal H}\mathcal G(X_{\varphi/\psi})=\dim_{\mathcal B}\mathcal
G(X_{\varphi/\psi})=d+1-\min(\bar\alpha_1^{-1},\dots,\bar\alpha_p^{-1})
\] in complete similarity to the result in \cite{ayache1}. Hence
we have constructed operator scaling Gaussian random fields with
the same box- and Hausdorff-dimension as the fractional Brownian
sheet, but our fields have additionally stationary increments.
\end{remark}

\noindent{\bf ACKNOWLEDGEMENT}\\[5pt]
We would like to thank David A. Benson for stimulating discussions
that inspired and focused the research presented in this paper.

\bibliographystyle{plain}

\end{document}